\documentclass{amsart}
\usepackage{amssymb,
enumitem,
mathrsfs,
hyperref,
tikz,
upgreek
}
\hypersetup{colorlinks=true}
\numberwithin{equation}{section}
\usepackage[T1]{fontenc}
\usepackage[shadow]{todonotes}

\newcommand{\RR}{\mathbb{R}}
\newcommand{\NN}{\mathbb{N}}
\newcommand{\QQ}{\mathbb{Q}}
\newcommand{\pre}[2]{#2^{#1}}

\newcommand{\Nbhd}{\boldsymbol{N}}

\newcommand{\markdef}[1]{\textbf{#1}}

\newcommand{\pow}{\mathscr{P}}
\newcommand{\Mod}{\operatorname{Mod}}

\newenvironment{enumerate-(a)}{\begin{enumerate}[label={\upshape (\alph*)}, leftmargin=2pc]}{\end{enumerate}}
\newenvironment{enumerate-(a)-r}{\begin{enumerate}[label={\upshape (\alph*)}, leftmargin=2pc,resume]}{\end{enumerate}}
\newenvironment{enumerate-(a)-5}{\begin{enumerate}[label={\upshape (\alph*)}, leftmargin=2pc,start=5]}{\end{enumerate}}
\newenvironment{enumerate-(A)}{\begin{enumerate}[label={\upshape (\Alph*)}, leftmargin=2pc]}{\end{enumerate}}
\newenvironment{enumerate-(A)-r}{\begin{enumerate}[label={\upshape (\Alph*)}, leftmargin=2pc,resume]}{\end{enumerate}}
\newenvironment{enumerate-(i)}{\begin{enumerate}[label={\upshape (\roman*)}, leftmargin=2pc]}{\end{enumerate}}
\newenvironment{enumerate-(i)-r}{\begin{enumerate}[label={\upshape (\roman*)}, leftmargin=2pc,resume]}{\end{enumerate}}
\newenvironment{enumerate-(I)}{\begin{enumerate}[label={\upshape (\Roman*)}, leftmargin=2pc]}{\end{enumerate}}
\newenvironment{enumerate-(I)-r}{\begin{enumerate}[label={\upshape (\Roman*)}, leftmargin=2pc,resume]}{\end{enumerate}}
\newenvironment{enumerate-(1)}{\begin{enumerate}[label={\upshape (\arabic*)}, leftmargin=2pc]}{\end{enumerate}}
\newenvironment{enumerate-(1)-r}{\begin{enumerate}[label={\upshape (\arabic*)}, leftmargin=2pc,resume]}{\end{enumerate}}
\newenvironment{itemizenew}{\begin{itemize}[leftmargin=2pc]}{\end{itemize}}

\newtheorem{theorem}{Theorem}[section]

\theoremstyle{definition}
\newtheorem{defin}[theorem]{Definition}
\newtheorem{example}[theorem]{Example}

\theoremstyle{remark}

\begin{document}

\title[Classification problems from the DST perspective]{Classification problems from \\ the descriptive set theoretical perspective}
\date{\today}
\author[L.~Motto Ros]{Luca Motto Ros}
\address{Dipartimento di matematica \guillemotleft{Giuseppe Peano}\guillemotright, Universit\`a di Torino, Via Carlo Alberto 10, 10121 Torino --- Italy}
\email{luca.mottoros@unito.it}
 \subjclass[2010]{Primary: 03E15}
 \keywords{Classification problems; complete invariants; Borel reducibility}
\thanks{The author would like to thank Filippo Calderoni and Matteo Viale for useful discussions and comments on a preliminary version of this paper. This work has been supported by the project PRIN~2017 ``Mathematical Logic: models, sets, computability'', prot. 2017NWTM8R}

\begin{abstract} 
Twenty years have passed since Kechris' seminal survey paper~\cite{Kechris1999}. As a follow-up of that work, we review some of the (anti-)classification results that have been obtained in the last decade using Borel reducibility and its generalizations to uncountable cardinals.
\end{abstract}

\maketitle

\section{Classification problems and Borel reducibility} \label{sec:intro}

One of the means by which human beings try to understand the world around them is classification. For example, we classify biological organisms in kingdoms, phyla, classes, etc., down to species and subspecies (taxonomy) to understand their characteristics and behaviors; and we classify elements in Mendeleev's periodic table to understand their chemical and physical properties. The classification process appears to be quite basic and to some extent innate,%
\footnote{Arguably, the natural numbers \( 1 \), \( 2 \), \( 3 \), \dots are abstracted through the process of classifying finite (small) sets up to one-to-one correspondence.}
and shows up also in daily life: from the ``classification'' of clothes in a wardrobe according to their characteristic to the classification of books in literary genres, and so on.

Mathematics is of course no exception: we e.g.\ classify connected closed surfaces (up to homeomorphism) using orientability and Euler characteristic,  Bernoulli automorphisms (up to conjugacy) by their entropy,  abelian \( p \)-groups (up to isomorphism) using Ulm invariants, and so on --- one may find beautiful and deep classification results in virtually every mathematical field. In contrast to the other scientific disciplines, however, mathematics has also the language and the tools to analyze its classification problems from a more abstract point of view. Indeed, any classification problem may be formalized as a set of objects  equipped with some natural equivalence relation  which aims at identifying objects sharing the same properties; a solution to such a classification problem is given by an assignment of complete invariants. More formally:

\begin{defin}
A \markdef{classification problem} is a pair \( (X,E) \) where \( X \) is a nonempty set and \( E \) is an equivalence relation on \( X \). A \textbf{solution} to a classification problem \( (X,E) \) is a function \( \varphi \) from \( X \) to some set \( I \) (whose elements will be called \textbf{invariants}) such that for all \( x,y \in X \)
\[ 
x \mathrel{E} y \iff \varphi(x) = \varphi(y).
 \] 
 The pair \( (\varphi,I) \) is called an \textbf{assignment of complete invariants} to the problem \( (X,E) \).
 \end{defin}
 
If the invariants in \( I \) are well-understood (e.g.\ they are natural or real numbers) and the function \(\varphi \) is easy to compute (e.g.\ it is continuous or Borel, if we are dealing with topological spaces), then such an assignment allows us to completely classify the elements of \( X \): given \( x,y \in X \), it is enough to compute \(\varphi \) and check whether \( \varphi(x) = \varphi(y) \) or not to determine if \( x \) and \( y \) are the ``same'' object (up to \( E \)).

\begin{example}
Fix \( n > 0 \). Let \( X \) be the set of all complex square matrices of order \( n \) and let \( E \) be the relation of similarity: a solution to this classification problem is given by canonical Jordan forms. More precisely, an assignment of complete invariants is the function \( \varphi \colon X \to I \) which computes for each element in \( X \) its Jordan form (so \( I \) is the set of all block diagonal \( n \)-matrices or, more precisely, of the sets of their blocks, as the Jordan form of a matrix is unique up to reordering of the Jordan blocks). 
\end{example}

This abstract formulation of classification problems and their solutions also suggests a way to assess their relative complexity. 

\begin{defin} \label{def:reduction}
Given two classification problems \( (X,E) \) and \( (Y,F) \),
we say that \( E \) \markdef{reduces} to \( F \) if there is a function \( f \colon X \to Y \), henceforth called a \markdef{reduction}, such that for all \( x,y \in X \)
\[ 
x \mathrel{E} y \iff f(x) \mathrel{F} f(y).
 \] 
 \end{defin}
Notice that, in particular, an assignment \( \varphi \colon X \to I \) of complete invariants to \( (X,E) \) may be viewed as a reduction of the latter to the identity relation on \( I \). If the function \( f \) in Definition~\ref{def:reduction} is again reasonably simple, then any solution of the problem \( (Y,F) \) may be turned into a solution of the problem \( (X,E) \) by composing the assignment of complete invariants to \( (Y,F) \) with the reduction \( f \). Of course, ``simplicity'' is a vague concept: if we want to analyze classification problems from such a general perspective, we must give a concrete meaning to this notion, and here is where descriptive set theory enters the picture.

In most cases, the sets of objects \( X \), \( Y \), \dots may be construed, up to coding, as completely metrizable second-countable (briefly, \markdef{Polish}) spaces, or even just \markdef{standard Borel} spaces. In these cases, it is then natural to consider the following notion.

\begin{defin}
Let \( E \) and \( F \) be equivalence relations on Polish/standard Borel spaces \( X \) and \( Y \), respectively. Then \( E \) is \markdef{Borel reducible} to \( F \) if there is a Borel function \( f \colon X \to Y \) reducing \( E \) to \( F \).
The preorder \( \leq_B \) is called \markdef{Borel reducibility}, its strict part is denoted by \( <_B \), while its induced equivalence relation is called \markdef{Borel bi-reducibility} and is denoted by \( \sim_B \).
\end{defin}

Borel reducibility was introduced in a seemingly independent way in \cite{Friedman1989} and \cite{HKL}. In the former, it was used as a tool to assess the complexity of first-order (or \( \mathcal{L}_{\omega_1 \omega} \)-)theories by looking at the topological complexity of the isomorphism relation on their countable models; in the latter, it was used to formulate and prove a generalization of a dichotomy by Glimm and Effros to all Borel equivalence relations.

Our aim is to survey some recent developments in the theory of Borel reducibility, especially in connection with various naturally occurring classification problems. The paper is organized as follows.
In Section~\ref{sec:borel} we discuss a first criterion (i.e.\ Borel vs.\ non-Borel) to distinguish tractable classification problems from the intractable ones, and provide concrete examples where this distinction has been employed.
In Section~\ref{sec:structural} we recall some of the more important dividing lines in the structure of (analytic) equivalence relations under Borel reducibility. This provides a rough taxonomy for classification problems and a way to describe their complexity. 
In Section~\ref{sec:examples} we  turn to some concrete classification problems that have been analyzed using this method. The list is very far from being complete, as Borel reducibility has been one of the most active fields in set theory in the last 30 years and produced an impressive amount of deep results. Nevertheless, we hope that our small selection will at least succeed in demonstrating the effectiveness of these techniques to tackle classification problems from the most diverse areas. 
Finally, in Section~\ref{sec:uncountable} we describe further generalizations of the notion of Borel reducibility that have been used to analyze classification problems for uncountable structures and non-separable spaces. This is a quite recent trend in the area, yet it already produced interesting and deep connections with e.g.\ model theory, algebra, and (functional) analysis.

Throughout the paper we assume familiarity with the basics of classical descriptive set theory as presented e.g.\ in~\cite{Kechris1995}. See also~\cite{gaobook}  for a comprehensive account on Borel reducibility.

\section{Borel vs.\ non-Borel} \label{sec:borel}

A first dividing line among classification problems \( (X,E) \) with \( X \) Polish (or standard Borel) is whether \( E \) is 
Borel as a subset of \( X \times X \) or not. Indeed, Borelness of \( E \) corresponds to having a 
procedure that, using only countable information/operations and in countably many steps, allows us to 
determine whether two given elements \( x,y \in X \) are \( E \)-equivalent. Of course this does not 
mean that we have a nice assignment of complete invariants for the problem at hand, but it at least 
shows that such a solution is not completely out of question. In contrast, knowing that \( E \) is 
\emph{not} Borel constitutes an evidence that the problem is very hard, and no simple solution can be 
expected; this kind of results should thus be regarded as anti-classification theorems. It should be noticed that the class of Borel classification problems is \( \leq_B \)-downward closed, thus it forms an initial fragment of the structure of all classification problems up to Borel reducibility.

A few results along this lines have been proposed already in~\cite{Friedman1989}. For example:
\begin{itemizenew}
\item
The following classification problems are Borel: 
countable finitely branching trees up to isomorphism,
countable finite-rank torsion-free abelian groups up to isomorphism,
countable fields of characteristic \( p \) and finite transcendence degree up to isomorphism (\cite[Theorems 2, 8 and 9]{Friedman1989}). 
\item
In contrast, for any prime \( p \) the  isomorphism relation on the class of countable abelian \( p \)-groups is complete analytic,%
\footnote{Recall that a subset \( A \) of a Polish space \( X \) is called \markdef{analytic} if it is a continuous image of a Polish space. All Borel sets are analytic, but there are analytic sets which are not Borel (indeed Borel sets are exactly the analytic sets with analytic complement). An analytic set is \markdef{complete} if all analytic sets are continuous preimages of it; thus, complete analytic sets are not Borel.}
 and thus not Borel (\cite[Theorem 6]{Friedman1989}).\label{thm:FSnonborel}
 \item
 The same applies to isomorphism on e.g.\ countable graphs, countable linear orders, and so on (see~\cite[Theorems 1, 3, 7 and 10]{Friedman1989}). However, in this cases one actually gets a more precise result (which  by~\cite[Theorem 4]{Friedman1989} implies the non-Borelness of the corresponding classification problem), namely, the \( S_\infty \)-completeness of such isomorphism relations --- see  Section~\ref{sec:structural} for the relevant definitions.
\end{itemizenew}

One of the most celebrated anti-classification results of this kind, due to Foreman, Rudolph and Weiss, is the following. (Readers not familiar with the mathematical objects involved may want to give a quick look at Section~\ref{subsec:measurepreserving} for some basic definitions before reading the statement below.)

\begin{theorem}[\cite{Foreman2011}] \label{thm:FRW}
The conjugacy relation on ergodic measure preserving transformations of
the unit interval with Lebesgue measure is complete analytic, and thus not Borel.
\end{theorem}

This result explains why the problem of determining whether
ergodic transformations are isomorphic or not has proven to be quite difficult, and is in contrast to the
situation of unitary operators where the spectral theorem can be used to
show that the conjugacy relation on the unitary group is Borel. In~\cite{Foreman2011} it is also shown that 
restricting the attention to rank one transformations, which form a generic
subset of the group of measure preserving transformations considered in the above theorem, is instead 
a Borel classification problem.

Similar anti-classification results have been obtained in the context of von Neumann factors (see also Section~\ref{subsubsec:vonNeumann} for more on this and in particular for the relevant definitions). 

\begin{theorem}[{\cite{Sasyk2009}}]
The isomorphism relation for von Neumann \( \mathrm{II}_1 \) factors is analytic but not Borel.
\end{theorem}

The following theorem deals instead with the automorphisms of a special \( C^* \)-algebra, that is, the Cuntz algebra \( \mathcal{O}_2 \) introduced in~\cite{Cuntz1977} (see Section~\ref{subsubsection:automorphisms} for more information on this and related problems).

\begin{theorem}[\cite{Gardella2016}] \label{thm:gard16}
The classification of the automorphisms of \( \mathcal{O}_2 \) up to  conjugacy and cocycle conjugacy is not Borel. The same conclusion holds even when one only considers automorphisms of \( \mathcal{O}_2 \) of finite order.
\end{theorem}

Theorems~\ref{thm:FRW}--\ref{thm:gard16}
are all establishing that certain classification problems are intractable.
In the opposite direction, we e.g.\ have the following (we refer the reader to Section~\ref{subsec:representations} for the relevant definitions).

\begin{theorem}[\cite{Hjorth2012}] 
The problem of classifying unitary representations of a given discrete countable group up to conjugacy is Borel. 
\end{theorem}

More generally, the same is true for representations of separable involutive Banach algebras, and so in particular it applies to separable \( C^* \)-algebras and to unitary representations of second countable locally compact groups.

\section{Structural results} \label{sec:structural}

The distinction between Borel and non-Borel classification problems is quite rough. A finer and more informative analysis can be obtained by considering how much complex a system of complete invariants must be in order to solve the classification problem at hand. For example, one of the simplest kind of complete invariants one may hope to find is constituted by real numbers or, equivalently, by the elements of any Polish or standard Borel space. (In what follows, \( \mathrm{id}(X) \) denotes the identity relation on \( X \).)

\begin{defin} \label{def:smooth}
A classification problem \( (X,E) \) is called \markdef{smooth} or \markdef{concretely classifiable} if it is Borel reducible to \( \mathrm{id}(\RR) \).
\end{defin}

This is the first important dividing line among (Borel) classification problems.
It turns out that there are only a few types of smooth equivalence relations (up to Borel bi-reducibility), and that such types are completely determined by the number of equivalence classes: finitely many, countably many, or perfectly many. This easily follows from the following dichotomy theorem, which is one of the first fundamental structural results in the theory.

\begin{theorem}[\cite{Silver1980}]
Let \( E \) be a Borel equivalence relation on a Polish space \( X \). Then either \( E \) has countably many classes (equivalently, \( E \leq_B \mathrm{id}(\NN) \)), or else \( \mathrm{id}(\RR) \leq_B E \).
\end{theorem}

(The above theorem actually applies to all coanalytic equivalence relations: however, we will not consider equivalence relations of that complexity here.)

A slightly more complicated kind of invariants is given by elements of a given Polish space up to ``countably many mistakes''. An example is given by reals up to rational translations, that is, by the Vitali equivalence relation \( E_v \) on \( \RR \) defined by
\[ 
x \mathrel{E_v} y \iff x-y \in \QQ.
 \] 
Up to Borel bi-reducibility, \( E_v \) is the same as the relation  \( E_0 \) of eventual equality on the Cantor space \( \pre{\omega}{2} \), where
\[ 
x \mathrel{E_0} y \iff \exists n  \, \forall m \geq n \, (x(m) = y(m)).
 \]

\begin{defin} \label{def:countableborel}
An equivalence relation on a Polish space is \markdef{countable Borel} if it is Borel and all its equivalence classes are (at most) countable.
A classification problem \( (X,E) \) is \markdef{essentially countable Borel} if it is Borel reducible to some countable Borel equivalence relation, that is, if it admits as complete invariants the equivalence classes of such a relation.
\end{defin}

Interestingly enough, among non-smooth Borel classification problems there is a minimal complexity (up to Borel bi-reducibility) corresponding exactly to \( E_0 \). This follows from the following dichotomy theorem by Harrington, Kechris and Louveau, which generalizes earlier work by Glimm and Effros.

\begin{theorem}[\cite{HKL}]
Let \( E \) be a Borel equivalence relation. Then either \( E \) is smooth, or else \( E_0 \leq_B E \).
\end{theorem}

The situation above \( E_0 \) is much wilder (see e.g.~\cite{Louveau1994,Adams2000}), and basically no other general dichotomy theorem as above can be obtained. In particular, Conley and Miller proved that  there is no ``immediate successor'' of \( E_0 \) even when restricting the attention to countable Borel equivalence relations, one of the simplest classes properly extending concretely classifiable problems.

\begin{theorem}[\cite{Conley2017}] 
Every basis for the countable Borel equivalence relations above \( E_0 \) is uncountable, that is: there is no countable family \( \{ F_n \mid n \in \mathbb{N} \} \) of equivalence relations such that for any countable Borel \( F \) with \( E_0 <_B F \) it must be the case that \( F_n \leq_B F \) for some \( n \in \mathbb{N} \).
\end{theorem}

The same is true even when weakening the notion of Borel reducibility to \emph{measure reducibility}, that is when the Borel reducibility is required to exist only on a co-null set (with respect to an appropriate measure).

An even more generous kind of invariants that has been considered in the literature, e.g.\ in connection to classification problems related to \( C^* \)-algebras (see Section~\ref{subsec:C*algebras}), is given by countable structures up to isomorphism, such as graphs, linear orders, trees, groups, and so on. (Recall that if \( \mathcal{L} \) is a countable relational language, then the countable \( \mathcal{L} \)-structures can be coded,%
\footnote{Basically, each structure on \(\omega\) is identified with the characteristic functions of its predicates.}
 up to isomorphism, as the elements of a suitable Polish space \( \mathrm{Mod}^\omega_\mathcal{L} \).)

\begin{defin} \label{def:classifablebycountablestructures}
A problem \( (X,E) \) is \markdef{classifiable by countable structures} if it is Borel reducible to the isomorphism relation on countable structures of a given first-order relational signature \( \mathcal{L} \).
\end{defin}

Notice that countable structures constitute a quite weak solution for a classification problem: indeed, 
among problems which are classifiable by countable structures there are some which are not even Borel 
(e.g.\ those mentioned on p.~\pageref{thm:FSnonborel}). Nevertheless, such a solution is considered 
perfectly acceptable when coming to very complicated classification problems, as illustrated by 
\( K \)-theoretic invariants for \( C^* \)-algebras --- see Section~\ref{subsec:C*algebras} 
for more on this. 

The largest class of classification problems usually analyzed through Borel reducibility is that of analytic equivalence relations (i.e.\ equivalence relations \( E \) which are analytic as subsets of the square \( X \times X \), where \( X \) is the Polish/standard Borel domain of \( E \)). As a matter of fact, the vast majority of classification problems naturally occurring in mathematics falls, up to a suitable coding procedure, inside this class. However, equivalence relations of this kind may be very complicated and intractable, so when possible it is customary to restrict the attention to the following strictly smaller class.

\begin{defin} \label{def:orbitrquivalencerelation}
An equivalence relation \( E \) on a Polish space \( X \) is an \markdef{orbit equivalence relation} if there is a Polish group \( G \) and a continuous (equivalently, Borel) action of \( G \) on \( X \) such that the \( E \)-equivalence classes are exactly the orbits of such action. A relation is \markdef{essentially orbit} if it is Borel reducible to an orbit equivalence relation.
\end{defin}

This already imposes some nontrivial complexity restrictions: for example, if \( E \) is an (essentially) orbit equivalence relation, then all its equivalence classes must be Borel. 

Since the isomorphism relation on \( \mathrm{Mod}^\omega_{\mathcal{L}} \) is induced by a continuous action of the symmetric group \( S_\infty \) on infinitely many elements,  i.e.\ the group of all permutations of \( \omega \), it turns out that all  kinds of classification problems considered so far (smooth, essentially countable Borel, classifiable by countable structures) are essentially orbit equivalence relations. Hjorth isolated a topological condition, called \emph{turbulence} (see~\cite{Hjorth2000}), to distinguish when an orbit equivalence relation is classifiable by countable structures or not: this is a key tool to show that a given classification problem does not admit countable structures (up to isomorphism) as invariants, thus yielding very strong anti-classification results.

Let us finally point out that all the classes of equivalence relations mentioned in this section admit complete elements, that is, equivalence relations  which are the most complicated ones in that class (up to Borel bi-reducibility).

\begin{defin} \label{def:complete}
Let \( \Gamma \) be a collection of equivalence relations. An equivalence relation \( E \) is \markdef{complete} (\markdef{for \( \Gamma \)}) if \( E \) is (Borel bi-reducible with an element) in \( \Gamma \) and \( F \leq_B E \) for all \( F \in \Gamma \).
\end{defin}

Here are some examples.
\begin{itemizenew}
\item
The orbit equivalence 
relation \( E_\infty \) induced by the shift action of the free group \( \mathbb{F}_2 \) with two generators on the space \( \pow(\mathbb{F}_2) \) of all subsets of \( \mathbb{F}_2 \) is complete for the class of (essentially) countable Borel equivalence relations.
\item
The relation \( \cong_{\mathsf{GRAPHS}} \) of isomorphism on countable graphs is \markdef{\( S_\infty \)-complete}, that is, it is complete for the class of problems which are classifiable by countable structures (the same applies e.g.\ to isomorphism on countable trees or countable linear orders).
\item
Let \( H([0;1]^\NN) \) be the automorphisms group of the Hilbert cube \( [0;1]^\NN \), and \( F(H([0;1]^{\NN}) \) be the associated Effros Borel space of closed subgroups of \( H([0;1]^\NN) \). Then
the orbit equivalence relation \( E^\infty_{H([0;1]^{\NN})} \) induced by the coordinatewise right action of \( H([0;1]^{\NN} )\) on the countable power of \( F(H([0;1]^{\NN}) \) is complete for the class of all (essentially) orbit equivalence relations. 
\end{itemizenew}

Results establishing that a given classification problem is complete for a certain class are nearly optimal, as they exactly determine which kind of invariants can (or cannot) be used to solve the problem. For example, if a problem \( (X,E) \) is \( S_\infty \)-complete, then it means that:
\begin{enumerate-(1)}
\item
such classification problem has a solution using countable structures (up to isomorphism) as invariants, but
\item
no strictly simpler objects may form a system of complete invariants for \( E \).
\end{enumerate-(1)}
We also mention that if a relation is \( S_\infty \)-complete or above, then it is necessarily proper analytic (and thus not Borel). 

Figure~\ref{fig:sumup} summarizes the dividing lines discussed so far, and may be used as a reference for the complexity of the examples considered in the next section.

\begin{center}
\begin{figure}
\begin{tikzpicture}[scale=1]
\draw [very thick] (0,0) to [out=10,in=270] (2.5,3.5) to [out=90,in=350] (0,7); 
\draw [very thick] (0,0) to [out=170,in=270] (-2.5,3.5) to [out=90, in=190] (0,7); 
\draw [fill] (0,0) circle [radius=0.07]; 
\draw [fill] (0,7) circle [radius=0.07]; 
\draw [red,dashed, ->] (2.5,6.5) to [out=150,in=5] (0.2,7.05);
\node [red, align=center,right] at (2.5,6.5) {complete analytic \\ equivalence relations};
\draw [red,dashed,<-] (0.2,-0.05) to [out=-20, in=-160] (1,-0.05);
\node [red,right] at (1,0) {\( \mathrm{id(\mathbb{N})} \)}; 
\draw [fill] (0,0.7) circle [radius=0.07]; 
\draw [red,dashed, <-] (0.2,0.75) to [out=10,in=155] (2,0.5);
\node [red,right] at (2,0.5) {\( \mathrm{id(\mathbb{R})} \)};
\draw (-1.02,0.4) to [out=32,in=185] (0,0.7) to [out=355,in=148] (1.02,0.4);
\draw [->] (-1.2,0.25) to (-0.4,0.35);
\node [left] at (-1.2,0.25) {smooth};
\draw [cyan, dotted, very thick] (-2.44,2.8) to [out=35,in=180]  (0,3.6) to [out=0,in=145] (2.44,2.8); 
\node [cyan] at (0,3.3) {Borel};
\node [cyan, align=center] at (0,4.1) {proper \\ analytic};
\draw (-1.51,0.8) to [out=80,in=185] (0,2.1) to [out=355,in=100] (1.51,0.8);
\draw [->] (-2.2,1.05) to (-0.8,1.3);
\node [align=center,left] at (-2.2,1.05) {essentially \\ countable Borel}; 
\draw [fill] (0,2.1) circle [radius=0.07]; 
\draw [red,dashed, <-] (0.2,2.15) to [out=10,in=160] (2.55,1.95);
\node [red,right] at (2.55,1.95) {\( E_{\infty} \)};
\draw [fill] (0,2.7) circle [radius=0.07]; 
\draw [red,dashed, <-] (0.2,2.75) to [out=10,in=160] (2.75,2.55);
\node [red,right] at (2.75,2.55) {\( \mathrm{id}(\RR)^+ \)};
\draw [fill] (0,1.4) circle [radius=0.07]; 
\draw [red,dashed, <-] (0.2,1.45) to [out=10,in=160] (2.35,1.25);
\node [red,right] at (2.35,1.25) {\( E_{0} \)};
\draw  (1.51,0.8) to [out=60,in=345]  (0,5); 
\draw  (-1.51,0.8) to [out=120,in=195] (0,5); 
\draw [fill] (0,5) circle [radius=0.07]; 
\draw [red,dashed, <-] (0.1,5.1) to [out=0,in=130] (3,3.7);
\node [red,right] at (3,3.7) {\( \cong_{\mathsf{GRAPHS}} \)};
\draw [->] (-2.8,2.5) to (-1,2.8);
\draw [->] (-2.8,2.5) to (-1,4.0);
\node [align=center, left] at (-2.8,2.5) {classifiable by \\ countable structures};
\draw  (1.51,0.8) to [out=50,in=353] (0,6); 
\draw  (-1.51,0.8) to [out=130,in=187] (0,6); 
\draw [fill] (0,6) circle [radius=0.07]; 
\draw [red,dashed, <-] (0.2,6.05) to [out=10,in=135] (2.8,5);
\node [red,align=center,right] at (2.8,5) {\( E^\infty_{ H([0;1]^\mathbb{N})} \)};
\draw [->] (-2.8,4.2) to (-1,5);
\node [align=center, left] at (-2.8,4.2) {essentially orbit \\ equivalence relations};
\draw [->] (-2.7,5.7) to (-0.8,6.2);
\node [align=center, left] at (-2.7,5.7) {other analytic \\ equivalence relations};
\end{tikzpicture}
\caption{The structure of analytic equivalence relations under Borel reducibility, together with some of the most important subclasses and some equivalence relations which are used as milestones when determining the complexity of classification problems.} \label{fig:sumup}
\end{figure}
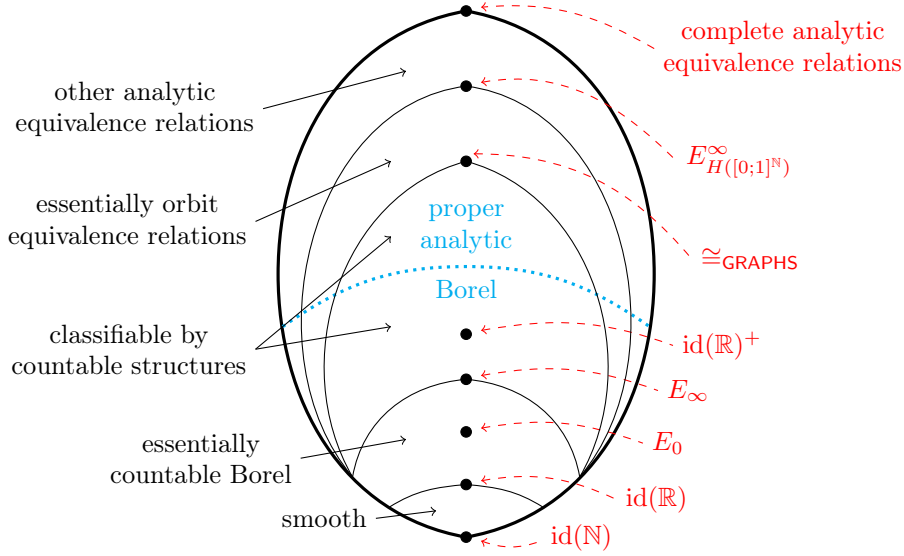
\end{center}

\section{Some examples} \label{sec:examples}

In this section we present  a (largely incomplete) selection of classification problems that have been analyzed using Borel reducibility. We mostly concentrate on results published in the last ten years --- for less recent results the reader may consult e.g.~\cite{Kechris1999} or~\cite{Foreman2018}.

\subsection{Groups}

Isomorphism on (various subclasses of) countable groups has been one of the mainstream lines of investigation since the very beginnings of the Borel reducibility theory. The literature on this subject is huge, so we are going to present just a sample of the results obtained.

Already in~\cite[Theorem 7]{Friedman1989} it was observed that, by work of Mekler, isomorphism on countable groups is \( S_\infty \)-complete (actually this is true when restricting to exponent \( p \) nilpotent class \( 2 \) groups, for any prime \( p \)). In contrast, by~\cite[Theorem 5]{Friedman1989} the isomorphism relation on abelian torsion groups is strictly \( \leq_B \)-simpler than that. The situation for torsion-free abelian groups is less clear, as only rank \( 1 \) groups have been classified in a satisfactory way. Hjorth showed that the isomorphism on rank \( 2 \) groups is strictly more complicated with respect to Borel reducibility, and then Thomas showed that the complexity of the corresponding classification problems \( \leq_B \)-strictly increases with the rank along all natural numbers. More recently, this has been refined by considering \( S \)-local groups (for \( S \) a set of primes), i.e.\ groups \( G \) such that \( p G = G \) for all primes \( p \notin S \).

\begin{theorem}[{\cite{Thomas2011}}]
Suppose that \( n \geq 2 \) and that \( S, T \) are sets of primes. Then the classification problem for the \( S \)-local torsion-free abelian groups of rank \( n \) is Borel reducible to the classification problem for the \( T \)-local torsion-free abelian groups of rank \( n \) if and only if \( S \subseteq T \).
\end{theorem}

Further results along the lines of the above theorem, together with information on the quasi-isomorphism relation, can be found in~\cite{Coskey2012}.

Finitely generated groups attracted a lot of attention as well. The isomorphism relation on such class is complete for countable Borel equivalence relations~\cite{Thomas1999v}. Recently, it has been observed that the latter relation is nontrivial already when restricting to complete groups, i.e.\ centerless groups with only inner automorphisms.

\begin{theorem}[{\cite{Thomas2019}}]
The isomorphism relation on finitely generated complete groups is not smooth.
\end{theorem}

As for bi-embeddability between groups we have:

\begin{theorem}[{\cite{Williams2014}}]
\begin{enumerate-(1)}
\item
The bi-embeddability relation on arbitrary countable groups is complete for all analytic equivalence relations (and thus it is strictly more complicated than isomorphism on the same class).
\item
The bi-embeddability relation on finitely generated groups is complete for countable Borel equivalence relations (and thus has the same complexity of isomorphism on the same class).
\end{enumerate-(1)}
\end{theorem}

The bi-embeddability relation on finitely generated groups has been further studied in~
\cite{Thomas2013,Thomas2016}. Calderoni and Thomas also carried out a systematic comparison between isomorphism and bi-embeddability on torsion and torsion-free abelian groups. In the torsion-free case we have:

\begin{theorem}[{\cite{Calderoni2019}}] \label{thm:TFAcountable}
The relation of bi-embeddability on countable torsion-free abelian groups is complete for all analytic equivalence relations (and thus it is strictly more complicated than isomorphism on the same class).
\end{theorem}

The situation for torsion abelian groups is more subtle.

\begin{theorem}[{\cite{Calderoni2019}}] \label{thm:calddelta12}
The relations of bi-embeddability and isomorphism on countable torsion groups are incomparable with respect to \( \leq_B \). However under mild large cardinal assumptions%
\footnote{The existence of a Ramsey cardinal is overkill.} 
and considering absolute \( \Delta^1_2 \)-reducibility \( \leq_{a \Delta^1_2} \) instead of Borel reducibility, the isomorphism relation on countable torsion groups becomes strictly more complex than the bi-embeddability relation on the same class.
\end{theorem}

A similar result holds for countable abelian \( p \)-groups for any prime \(p \). The reducibility \( \leq_{a \Delta^1_2} \) is a technical yet quite natural weakening of \( \leq_B \) which sometimes gives a picture which better reflects the intuitively correct results about classification --- an illuminating example of this is Theorem~\ref{thm:calddelta12} itself. Using such reducibility, Shelah and Ulrich recently made a substantial progress on the major open problem of establishing the exact complexity of isomorphism on abelian torsion-free groups  (of arbitrarily high rank).

\begin{theorem}[{\cite{Shelah2019}}] \label{thm:shelahonTFAG}
It is consistent with \( \mathsf{ZFC} \) that the isomorphism on countable torsion-free abelian groups is \( S_\infty \)-complete with respect to the reducibility \( \leq_{a \Delta^1_2} \).
\end{theorem}

While this paper was in preparation, it has been announced in~\cite{paolini2021torsionfree} that isomorphism on countable torsion-free abelian groups is indeed \( S_\infty \)-complete (with respect to standard Borel reducibility). This result would thus supersede Theorem~\ref{thm:shelahonTFAG} and close a problem which resisted several attempts of solution since its first formulation in~\cite{Friedman1989}.

\subsection{Complete theories and model theory}

Given a first-order theory \( T \), we denote by \( \cong^\omega_T \) the isomorphism relation on (the codes of) its countable models, i.e.\ the restriction of \( \cong \) to \( \Mod^\omega_T = \{ A \in \Mod^\omega_{\mathcal{L}} \mid A \models T \} \) (where \( \mathcal{L} \) is the signature of \( T \)). Since \( \Mod^\omega_T \) is a Borel subset of the Polish space \( \Mod^\omega_{\mathcal{L}} \), the isomorphism relation \( \cong^\omega_T \) may be construed as an analytic equivalence relation.

\begin{theorem}[{\cite{Coskey2010}}]
Let \( T \) be any countable complete first-order theory extending Peano Arithmetic \( PA \). Then \( \cong^\omega_T \) is \( S_\infty \)-complete.
\end{theorem}

In~\cite{Laskowski2007}, Laskowski established the first link between stability theory and Borel reducibility. In particular he considered \(\omega\)-stable countable complete first-order theories and isolated a model-theoretic dividing line distinguishing between complicated (ENI-DOP) and more tractable (ENI-NDOP) theories. Then he proved that ENI-DOP theories \( T \) are complicated also from the descriptive set-theoretic point of view, namely, that \( \cong^\omega_T \) is \( S_\infty \)-complete for any such theory \( T \).

Then in~\cite{Koerwien2009} and, independently, in~\cite{Laskowski2015} it was defined a natural ``decomposition'' rank, called eni-depth, for models of ENI-NDOP theories. Such rank is either a countable ordinal, in which case the theory is called eni-shallow, or it is undefined, in which case the theory is called eni-deep. Combining the main results of the two papers one gets a quite clear picture remarkably connecting the eni-depth of a theory to the complexity of \( \cong^\omega_T \) with respect to Borel reducibility.

\begin{theorem} \label{thm:koerwien}
Let \( T \) be an \(\omega\)-stable ENI-NDOP countable complete first-order theory.
\begin{enumerate-(1)}
\item \label{thm:koerwien-1}
It \( T \) is eni-deep, then \( \cong^\omega_T \) is \( S_\infty \)-complete (\cite{Laskowski2015}).
\item \label{thm:koerwien-2}
If \( T \) is eni-shallow with eni-depth \( 1 \), then \( \cong^\omega_T \) is smooth (\cite{Koerwien2009}).
\end{enumerate-(1)}
Moreover 
\begin{enumerate}[label=\upshape(3), leftmargin=2pc]
\item \label{thm:koerwien-3}
There is an \( \omega_1 \)-chain of ENI-NDOP eni-shallow theories \( T_\alpha \) with strictly increasing eni-depth  such that \( (\cong^\omega_{T_\alpha})_{\alpha < \omega_1} \) form a \( \leq_B \)-increasing chain of equivalence relations cofinal among Borel isomorphism relations. (\cite{Koerwien2009}).
\end{enumerate}
\end{theorem}

By previous work of Shelah, Harrington and Makkai~\cite{Shelah1984}, we also know that if the eni-depth of an ENI-NDOP theory \( T \) is at least \( 3 \), then there are \( 2^{\aleph_0} \)-many non-isomorphic countable models of \( T \) (and the bound of \( 3 \) is optimal). In the same work, it is also shown that the Vaught conjecture holds for \( \omega \)-stable theories \( T \), i.e.\ that the number of non-isomorphic countable models of \( T \) is either \( \leq \aleph_0 \) or \( 2^{\aleph_0} \).

Despite the initial success depicted above of this line of research, Koerwien came up with an example of a very simple \(\omega\)-stable theory \( T \) for which \( \cong^\omega_T \) is not simple at all. This shattered hopes of finding a tight correspondence between simplicity from the the point of view of model theory and simplicity from the point of view of descriptive set theory, at least for what concerns countable models.

\begin{theorem}[{\cite{Koerwien2011}}]
There is an \(\omega\)-stable ENI-NDOP countable complete first-order theory \( T \) with eni-depth \( 2 \)  such that \( \cong^\omega_T \) is not Borel.
\end{theorem}

Koerwien's example has been modified in~\cite{Ulrich2017} to get a theory \( T \) with the same model-theoretic properties but such that \( \cong^\omega_T \) is furthermore \( S_\infty \)-complete. In the same paper, the authors also develop new techniques to provide various examples of complete theories \( T \) for which \( \cong^\omega_T \) is proper analytic yet not \( S_\infty \)-complete.

The correspondence between stability theoretic notions and the descriptive set-theoretic complexity of the isomorphism relation has been recently reconsidered in the context of uncountable models --- see Section~\ref{sec:uncountable}.

In a different direction, Rast analyzed  the complexity of isomorphism on countable models of the complete theories of suitable expansions of linear orders.

\begin{theorem}[{\cite{Rast2017a}}]
Let \( A \) be an expansion of a linear order by finitely many unary predicates. Then either \( \cong^\omega_{\mathrm{Th}(A)} \) has just one class (i.e.\ \( \mathrm{Th}(A) \) is \( \aleph_0 \)-categorical) or it is \( S_\infty \)-complete.
\end{theorem}

In particular, the above theorem applies to linear orders (with no other unary predicate) and to colored linear orders using finitely many colors. The paper~\cite{Rast2017a} analyzes also the case of countably many unary predicates added to a linear order (equivalently, colored linear orders using infinitely many colors), showing that there are just five possibilities for the complexity of \( \cong^\omega_{\mathrm{Th}(A)} \) and characterizing  in simple model-theoretic terms
when each case can occur.

Rast and Singh Sahota classified instead the possible complexities of \( \cong^\omega_T \) when \( T \) is a complete \( o \)-minimal theory, depending on its model-theoretic properties.

\begin{theorem}[{\cite{Rast2017}}]
Let \( T \) be a complete \( o \)-minimal theory in a countable language.
\begin{enumerate-(1)}
\item
If \( T \) has no non-simple types and the set \( S_1(T) \) (of types of \( T \) with at most one free variable) is countable, then \( \cong^\omega_T \) is smooth.
\item
If \( T \) has no non-simple types and the set \( S_1(T) \) is uncountable, then \( \cong^\omega_T \) is Borel bi-reducible with \( \mathrm{id}(\RR)^+ \)  (see Figure~\ref{fig:sumup}).%
\footnote{The symbol \( \mathrm{id}(\RR)^+ \) denotes the relation of equality between countable sets of reals: it is Borel, but strictly \( \leq_B \)-above all countable Borel equivalence relations.}
\item
If \( T \) admits a non-simple type, then \( \cong^\omega_T \) is \( S_\infty \)-complete.
\end{enumerate-(1)}
\end{theorem}

\subsection{Irreducible representations of (countable) groups} \label{subsec:representations}

Let \( U(\mathcal{H}_\infty) \) be the unitary group of isomorphisms of a separable infinite-dimensional Hilbert space \( \mathcal{H}_\infty \). Given a countable group \( G \), let \( \mathrm{Irr}_\infty(G) \) be the Polish space of all infinite-dimensional irreducible unitary representations of \( G \), that is, the space of all homomorphisms \( \rho \colon G \to U(\mathcal{H}_\infty) \) such that there is no nontrivial Hilbert subspace of \( \mathcal{H}_\infty \) which is invariant under all the \( \rho(g) \) for \( g \in G \). The group \( U(\mathcal{H}_\infty) \) acts by pointwise conjugacy on \( \mathrm{Irr}_\infty(G) \): given \( T \in U(\mathcal{H}_\infty) \) and \( \rho \in \mathrm{Irr}_\infty(G) \), the representation \( T \cdot \rho \in \mathrm{Irr}_\infty(G) \) is defined by
\[ 
(T \cdot \rho)(g) = T \circ \rho(g) \circ T^{-1}.
 \] 
We say that  \( \rho_1, \rho_2 \in \mathrm{Irr}_\infty(G) \) are \markdef{equivalent}, and write \( \rho_1 \approx_G \rho_2 \), if \( \rho_2 = T \cdot \rho_1 \) for some \( T \in U(\mathcal{H}_\infty) \). The \markdef{unitary dual \( \widehat{G} \)} of \( G \) is the quotient \( \mathrm{Irr}_\infty(G) /	\approx_G \) equipped with its Mackey Borel structure. It can be shown using~\cite[Theorem 2.5]{MR12} that if \( G \) and \( H \) are countable groups, then \( {\approx_G} \sim_B {\approx_H} \) if and only if the unitary duals \( \widehat{G} \) and \( \widehat{H} \) are Borel isomorphic.

The problem of classifying infinite-dimensional irreducible unitary representations of a countable group up to equivalence dates back to the Sixties. For example, by work of Glimm and Thoma the relation \( \approx_G \) is not smooth whenever \( G \) is not abelian-by-finite (equivalently, \( \mathrm{Irr}_\infty(G) \neq \emptyset \), or \( G \) is a non-type I group). Later on Hjorth  showed that for such \( G \) the relation \( \approx_{G} \) is not even classifiable by countable structures (this was one of the first applications of his turbulence theory, see~\cite[Chapter 5]{Hjorth2000} for a proof in the special case where \( G = \mathbb{F}_\infty \) is the free group with countably many generators). 

More recently, Thomas made further progress on this kind of problems. Following~\cite[Definition 1.8]{Thomas2015}, say that a countable group \( G \) is \markdef{representation universal} if \( {\approx_H} \leq_B {\approx_G } \) for all countable groups \( H \). (In other words, the unitary dual \( \widehat{G} \) is as complicated as possible.) It is not hard to see that \( \mathbb{F}_\infty \) is an example of a representation universal group. A less trivial example is the following:

\begin{theorem}[{\cite{Thomas2015}}]
The free group on two generators \( \mathbb{F}_2 \) is representation universal.
\end{theorem}

In the same paper it is also shown that amenable groups are \( \leq_B \)-minimal among non-type I groups.

\begin{theorem}[{\cite{Thomas2015}}]
If \( G \) is a countable non-type I group and \( H \) is a countable amenable
group, then \( {\approx_H} \leq_B {\approx_G} \). In particular, if \( G \) and \( H \) are countable
amenable non-type I groups, then \( \approx_G \) and \( \approx_H \) are Borel bi-reducible.
\end{theorem}

It is not yet known whether amenable non-type I groups are representation universal, although some partial results in~\cite{Thomas2015} suggest that this should not be the case.

\subsection{Metric spaces}

The classification of separable complete (briefly: Polish) metric spaces up to isometry dates back to the work of Gromov, who proved that such a classification problem is smooth when restricted to compact spaces. However, the general problem turned out to be much more complicated. For example, by work of Clemens, Gao and Kechris  the isometry relation on all Polish metric spaces is complete for orbit equivalence relations~\cite{Clemens2001, Gao2003}. Several natural subclasses of Polish metric spaces have been considered as well, including ultrametric spaces,%
\footnote{An ultrametric is a metric \(d \) satisfying the following strengthening of the triangular inequality: \( d(x,y) \leq \max \{ d(x,z), d(z,y) \} \) for all \( x,y,z \) in the space.} 
zero-dimensional spaces, discrete spaces, locally compact spaces, and so on.

Among the problems left open by~\cite{Gao2003} there was that of establishing  whether zero-dimensional Polish metric spaces can be classified using countable structures as invariants. This was answered negatively by Clemens.

\begin{theorem}[{\cite{clemensisometry}}]
Zero-dimensional Polish metric spaces are not classifiable up to isometry by countable structures.
\end{theorem}

In~\cite{clemensisometry} it is also shown that the same is true e.g.\ for ultrahomogeneous Polish metric spaces. (A metric space is ultrahomogeneous if any partial isometry between finite subsets of it can be extended to
an isometry of the whole space.)

Another problem left open by~\cite{Gao2003} was that of classifying locally compact Polish ultrametric spaces up to isometry. These are certainly classifiable by countable structures, as all Polish ultrametric spaces are, but there was some evidence that the complexity of the problem could be simpler than that. It turns out that this is not the case.

\begin{theorem}[{\cite{Camerlo2018}}] \label{thm:loccompactultrametric}
The relation of isometry between locally compact Polish ultrametric spaces is \( S_\infty \)-complete. In fact, the same is true even when restricting to discrete Polish ultrametric spaces.
\end{theorem}

In~\cite{louros}, Louveau and Rosendal also considered the relation of isometric bi-embeddability on Polish metric spaces. In particular, they showed that such relation is complete for all analytic equivalence relations even when restricted to discrete spaces or to ultrametric ones. However, whether the same applies to the class of spaces which are \emph{simultaneously} discrete (or even just locally compact) and ultrametric remained open for a while. Exploiting the techniques developed for Theorem~\ref{thm:loccompactultrametric} one gets: 

\begin{theorem}[{\cite{Camerlo2018}}] \label{thm:loccompactultrametric-2}
The relation of isometric bi-embeddability between discrete Polish ultrametric spaces is complete for all analytic equivalence relations. Thus the same is true for locally compact Polish ultrametric spaces.
\end{theorem}

An interesting feature of the analysis leading to Theorems~\ref{thm:loccompactultrametric}--\ref{thm:loccompactultrametric-2} is that the complexity of the isometry and isometric bi-embeddability relations on Polish ultrametric spaces depends only on the order type of the distances used in the space at hand: this is heavily in contrast with the situation for general Polish spaces --- see e.g.~\cite{Camerlo2019}.

Of course other equivalence relations can be considered on the class of Polish metric spaces: homeomorphism, uniform homeomorphism, bi-Lipschitz homeomorphism and so on. In this direction we have, for example:

\begin{theorem}[{\cite{Ferenczi:2009fk}}]
The relation of uniform homeomorphism on Polish metric spaces is complete for all analytic equivalence relations.
\end{theorem}

\subsection{Topological spaces and geometry}

The classification of geometrical objects is certainly among the leading problems in the area. It is no surprise that significant results were obtained already in the first years after the introduction of the notion of Borel reducibility. For example, Hjorth and Kechris~\cite{Hjorth2000bis} showed that 
the classification problem for Riemann surfaces considered up
to conformal (i.e.\ biholomorphic) equivalence is complete for countable Borel equivalence relations, while restricting the attention to \emph{compact} Riemann surfaces one gets a smooth equivalence relation. Moreover, for \( n \geq  2 \), the \( n \)-dimensional complex manifolds are not classifiable (up to biholomorphic equivalence) by countable structures.
 (To the best of our knowledge, it is still open whether these equivalence relations are Borel or not.)

\subsubsection{Knots}

A more recent result concerns the classification of knots. Recall that a \markdef{knot} is an embedding of a circle \( S^1 \) into the three-dimensional Euclidean space \( \RR^3 \). The natural equivalence between knots is the one induced by ambient isotopies or, equivalently, by the existence of an orientation-preserving homeomorphism from \( \RR^3 \) into itself sending a knot into the other one. An important dividing line among knots is the one distinguishing \markdef{tame} knots (i.e.\ knots equivalent to some polygonal) from \markdef{wild} knots. From the Borel reducibility point of view, the classification of tame knots  is trivial because there are only countably many of them (up to equivalence). In contrast:

\begin{theorem}[{\cite{Kulikov2017}}]
The equivalence relation on wild knots is not classifiable by countable structures, and it is indeed \( \leq_B \)-strictly above the \( S_\infty \)-complete equivalence relation.
\end{theorem}

This result mathematically explains why there is no hope to obtain a satisfactory classification for wild knots.

\subsubsection{Compact metric spaces}

Consider the problem of classifying \markdef{compact metric spaces} up to homeomorphism (rather than isometry). The relevance of the problem partially relies on the fact that it constitutes a lower bound for the problem of classifying separable \( C^* \)-algebras.
Hjorth~\cite{Hjorth2000} showed that such problem is not classifiable by countable structures, and it is indeed \( \leq_B \)-strictly above the \( S_\infty \)-complete equivalence relation. However, the exact complexity of the problem was settled only a few years ago.

\begin{theorem}[{\cite{Zielinski2016}}] \label{thm:zielinski}
The relation of homeomorphism on compact metric spaces is complete for orbit equivalence relations. \end{theorem}

As mentioned, this has direct consequences on the classification of (certain classes) of separable \( C^* \)-algebras (see Theorem~\ref{thm:zielinskiC^*algebras}). Theorem~\ref{thm:zielinski} has been generalized in~\cite{Rosendal2018}, where it is shown that the relation of homeomorphic isomorphism between compact metrizable \( \mathcal{L} \)-structures (where \( \mathcal{L}\) is any first-order countable relational language) is complete for orbit equivalence relations. This led to interesting information on the classification of topological groups --- see Theorem~\ref{thm:rosendalzielinski}.

\subsubsection{Continua}

A somewhat related problem is that of classifying \markdef{continua} (i.e.\ nonempty compact connected metric spaces) up to homeomorphism. They are all embeddable in the Hilbert cube \( [0;1]^\NN \), and if we restrict to continua of dimension \( n \in \NN \) then we can construe them as closed subsets of \( [0;1]^n \). It is easy to see that if \( n = 1 \) then there are only two classes of one-dimensional continua up to homeomorphism. In contrast,  it was proved in~\cite{Camerlo2005} that graph isomorphism is Borel reducible to homeomorphism of continua of dimension \( n \geq 2 \) (even when restricting to certain very simple continua called dendrites). This has recently been improved as follows.

\begin{theorem}[{\cite{Chang2019}}]
Continua of a fixed dimension \( n \geq 2 \) are not classifiable (up to homeomorphism) by countable structures. Indeed, the homeomorphism relation on such class is \( \leq_B \)-strictly above the \( S_\infty \)-complete equivalence relations.
\end{theorem}

This implies that the same is true when considering arbitrary closed subsets of \( [0;1]^n \). On such class one may also consider a different equivalence relation, namely the one induced by ambient homeomorphisms (i.e.\ homeomorphisms from \( [0;1]^n \) to itself sending a closed set to the other one) --- this is clearly reminiscent of the natural notion of equivalence for knots considered above.

\begin{theorem}[{\cite{Chang2019}}]
Closed sets of a fixed dimension \( n \geq 2 \) are not classifiable by countable structures up to ambient homeomorphism (such relation is again \( \leq_B \)-strictly above the \( S_\infty \)-complete equivalence relations).
\end{theorem}

Notice that for dimension \( 1 \), both the homeomorphism relation and the ambient homeomorphism relation on closed subsets of \( [0;1] \) are \( S_\infty \)-complete, and thus a bit simpler than the analogous relations for higher dimensions.

\subsection{Polish groups}

A topological group is \markdef{Polish} when its topology is separable and completely metrizable. The natural identification for such groups is topological isomorphism. The following result settles the classification problem for arbitrary Polish groups.

\begin{theorem}[{\cite{Ferenczi:2009fk}}]
The relation of topological isomorphism on (abelian) Polish groups is complete for all analytic equivalence relations.
\end{theorem}

Restricting the attention to certain natural classes of Polish groups, one gets instead nontrivial upper bounds. Recall that a topological group \( G \) is \markdef{Roelcke precompact} if for each open neighborhood \( U \) of the identity there is a finite set \( F \) such that \( UFU=G\).

\begin{theorem}[{\cite{Rosendal2018}}] \label{thm:rosendalzielinski}
The relation of topological isomorphism between locally compact and Roelcke precompact Polish groups is classifiable by compact metrizable structures, and thus these are essentially orbit equivalence relations.
\end{theorem}

In the same vein, Kechris, Nies and Tent considered the problem of classifying closed subgroups of \( S_\infty \) (endowed with the topology inherited from the Baire space) or, equivalently, the non-archimedean
Polish groups.
Besides the cases of locally compact%
\footnote{Notice that the class of locally compact closed subgroups of \( S_\infty \) represents, up to homeomorphism, the class of all totally disconnected locally compact Polish groups.}
 and Roelcke precompact groups, in this context it is natural to consider also
\markdef{oligomorphic} groups, i.e.\ closed subgroups \( G \) of \( S_\infty \) such that for each \( n \) the canonical action of \( G \) on the cartesian product \( \omega^n \) has only finitely many orbits.

\begin{theorem}[{\cite{Kechris2018}}]
The topological isomorphism relation for compact closed subgroups of \( S_\infty \) is \( S_\infty \)-complete. The same is true when considering locally compact or Roelcke precompact closed subgroups of \( S_\infty \).

Oligomorphic closed subgroups of \( S_\infty \) are classifiable by countable structures up to topological isomorphism.
\end{theorem}

\noindent
(The fact that the relation of isomorphism between locally compact or Roelcke precompact non-archimedean Polish groups is Borel reducible to isomorphism between countable structures was independently proved in~\cite{Rosendal2018}.)

As noticed in~\cite{Kechris2018}, the above classes of groups are naturally linked to notions in algebra or model theory. For example, it is
well-known that the closed subgroups of \( S_\infty \) are, up to topological group isomorphism, the automorphism groups
of countable structures,
while the oligomorphic groups are
precisely the automorphism groups of \(\omega\)-categorical structures with domain the
natural numbers. Notably, under the latter correspondence topological isomorphism turns into
bi-interpretability of the structures by a  result in Ahlbrandt and Ziegler~\cite{AZ84} going
back to unpublished work of Coquand.

As for topological bi-embeddability between Polish groups we instead have:

\begin{theorem}[{\cite{CMR18}}]
The relation of topological bi-embeddability on Polish groups is complete for all analytic equivalence relations. The same is true if we restrict the attention to Polish groups carrying a bi-invariant bounded metric under the relation induced by metric preserving group embeddings.
\end{theorem}

Thus the problems of classifying Polish groups up to topological isomorphism or up to topological bi-embeddability have the same complexity --- they are both completely intractable.

\subsection{Banach spaces}

The main notions of identification for \markdef{Banach spaces} (i.e.\ complete normed vector spaces) are linear isometry and isomorphism.
Not surprisingly, both relations are quite complicated, although when restricting to separable Banach spaces they have strictly different complexities.

\begin{theorem}[{\cite{Melleray2007}}]
The relation of (linear) isometry between separable Banach spaces is complete for orbit equivalence relations.
\end{theorem}

\begin{theorem}[{\cite{Ferenczi:2009fk}}]
The relation of isomorphism on separable Banach spaces is complete for all analytic equivalence relations. 
\end{theorem}

The same conclusion holds for the relations of Lipschitz isomorphism and complemented bi-embeddability.
As for bi-embeddability between separable Banach spaces, Louveau and Rosendal proved already in~\cite{louros} that the corresponding classification problem has maximal complexity: indeed,
the relation of linear isometric bi-embeddability between separable Banach spaces is complete for all analytic equivalence relations.

\subsection{Bounded operators and \( C^* \)-algebras} \label{subsec:C*algebras}

An area where some of the most spectacular (anti-)classification results have been obtained through Borel reducibility is the one concerning \markdef{\( C^* \)-algebras}, i.e.\ Banach algebras together with an involution map \( x \mapsto x^* \) satisfying the properties of the adjoint.

\subsubsection{Classification of bounded operators}

A first theme in the area is that of classifying \markdef{bounded operators} on an infinite-dimensional
separable complex Hilbert space \( \mathcal{H}_\infty \) up to various notions of equivalence. For example, Weyl and von Neumann proved that two bounded self-adjoint operators are unitarily equivalent modulo compact (i.e.\ one differs from a unitary conjugate of the other one by a compact operator) if and only if they share the same essential spectrum. As shown in~\cite{Ando2015}, this can be recast in terms of Borel reducibility by saying that the classification of bounded self-adjoint operators up to unitary equivalence modulo compact is smooth. In contrast

\begin{theorem}[{\cite{Ando2015}}] 
The  unbounded self-adjoint operators are not classifiable by countable structures up to unitary equivalence modulo compact.
\end{theorem}

Previously, Kechris and Sofronidis~\cite{Kechris2001} had shown that the same conclusion holds when classifying self-adjoint or unitary bounded operators up to unitary equivalence (rather than unitary equivalence modulo compact). Such operators are also not classifiable up to unitary equivalence using Borel actions of CLI groups, i.e.\ Polish groups admitting a compatible complete left-invariant metric (this includes all locally compact and all solvable Polish groups, hence in particular all abelian Polish groups).

\begin{theorem}[{\cite{Lupini2018g}}]
The relation of unitary equivalence on self-adjoint (alternatively, unitary) bounded operators is not Borel reducible to a Borel action of a CLI group on a Polish space.
\end{theorem}

Recently, Smythe considered other notions of equivalence for bounded operators, all defined by declaring two operators equivalent if their difference is in a certain prescribed class.

\begin{theorem}[{\cite{Smythe2017}}] 
\begin{enumerate-(1)}
\item
Equivalence modulo finite rank operators  is a Borel equivalence relation that is not essentially orbit (even when restricted to compact operators).
\item
 Equivalence modulo compact operators is a
Borel equivalence relation that is not classifiable by countable structures.
\item
Equivalence modulo Schatten \( p \)-class is
a Borel equivalence relation that is not classifiable by countable structures (even when restricted to compact operators).
\end{enumerate-(1)}
\end{theorem}

\subsubsection{Von Neumann factors} \label{subsubsec:vonNeumann}
A \markdef{von Neumann algebra} is a special type of \( C^* \)-algebra, namely, a \( * \)-algebra of bounded operators on a Hilbert space that is closed in the weak operator topology and contains the identity operator. 
A von Neumann algebra whose center is trivial (i.e.\ it consists only of multiples of the identity operator) is called a \markdef{factor}. Von Neumann showed in 1949 that every von Neumann algebra on a separable Hilbert space is isomorphic to a direct integral of factors, and such decomposition is essentially unique. Thus the problem of classifying isomorphism classes of von Neumann algebras on separable Hilbert spaces can be reduced to that of classifying isomorphism classes of factors.

Factors are usually divided into three classes (each of which has then several interesting subclasses):
\begin{itemizenew}
\item
\emph{Type I factors.} These are the factors that have a minimal projection \( E \neq 0 \), i.e.\ a projection \( E \) such that there is no other projection \( F \) with \( 0 < F < E \). Any factor of type I is isomorphic to the von Neumann algebra of all bounded operators on some Hilbert space: since there is only one Hilbert space for every cardinal number, this shows that it is quite easy to classify type I factors.
\item
\emph{Type II factors}. A factor is said to be of type II if there are no minimal projections but there are non-zero finite projections. This class is divided into the subclasses \( \mathrm{II}_1 \) and \( \mathrm{II}_\infty \), depending on whether the identity operator of the factor is finite or not. Factors of type \( \mathrm{II}_\infty \) can also be characterized as the tensor products of a factor of type \( \mathrm{II}_1 \) and an infinite type I factor.
\item
\emph{Type III factors.} These are the factors that do not contain any nonzero finite projection at all (their existence was proved by von Neumann in 1940). The class is then split into infinitely many subclasses \( \mathrm{III}_\lambda \), where \( 0 \leq \lambda \leq 1 \).
\end{itemizenew}

Woods showed in 1973 that the classification of separable type III factors is not smooth. Exploiting Borel reducibility and deep results of Popa, Sasyk and T\"ornquist greatly improved this by showing that indeed all non-type I factors are essentially unclassifiable.

\begin{theorem}[{\cite{Sasyk2009, Sasyk2010,Sasyk2019}}]
The separable factors in each of the classes \( \mathrm{II}_1 \), \( \mathrm{II}_\infty \), \( \mathrm{III}_\lambda \) (\( 0 \leq \lambda \leq 1 \)) are not classifiable by countable structures. The same applies even to the class of (free) Araki-Woods factors, and thus to the class of injective \( \mathrm{III}_0 \) factors. 
\end{theorem}

Recently, Spaas provided also interesting results concerning the classification of the Cartan subalgebras of a given von Neumann factor.

\begin{theorem}[{\cite{Spaas2018}}]
There is a large family of (non-hyperfinite) \( \mathrm{II}_1\) von Neumann factors whose Cartan subalgebras
 are not classifiable by countable
structures up to unitary conjugacy.
\end{theorem}

Spaas also constructs examples of \( \mathrm{II}_1 \) von Neumann factors whose Cartan subalgebras
up to conjugacy by an automorphism are not classifiable by countable structures, directly shows that the Cartan subalgebras of the hyperfinite \( \mathrm{II}_1 \) von Neumann factor up to unitary
conjugacy are not classifiable by countable structures, and
deduces that the same holds for any McDuff  \( \mathrm{II}_1 \) factor with at
least one Cartan subalgebra.

\subsubsection{Separable \( C^* \)-algebras}

An \markdef{approximately finite-dimensional} (briefly, \markdef{AF}) \markdef{\( C^* \)-algebra} is a \( C^* \)-algebra which is a direct limit \( \displaystyle{A = \lim_{\longrightarrow} (A_n , \varphi_n)} \) of an inductive system consisting of finite-dimensional \( C^* \)-algebras \( A_n \) and \( * \)-homomorphisms \( \upvarphi_n \) between them. Elliot gave a satisfactory classification of AF algebras developing the so-called \( K \)-theory:  Two AF algebras \( A \) and \( B \) are isomorphic if and only if their dimension groups \( (K_0(A), K_0^+(A), \Gamma(A)) \) and \( (K_0(B), K_0^+(B), \Gamma(B)) \) are isomorphic. (It may be worth recalling here that \( K_0(A) \) is an abelian ordered group whose positive cone is \( K_0^+(A) \), while  \( \Gamma(A) \) is the so-called scale of \( K_0(A) \), i.e.\ the collection of elements represented by projections in \( A \).) 
It turns out that Elliot's assignment of complete invariants can be computed in a Borel manner, once all the relevant objects are coded in suitable standard Borel spaces (which in this case is a highly nontrivial observation). It follows that one can recast Elliot's classification theorem in the Borel reducibility framework as

\begin{theorem}[{\cite{Farah2013}}]
AF \( C^* \)-algebras are classifiable by countable structures.
\end{theorem}

Another class of \( C^* \)-algebras for which a classification theory by means of \( K \)-theory and traces was (partially) developed is that of simple \markdef{nuclear \( C^* \)-algebras}, where a  \( C^* \)-algebra \( A \) is nuclear if the injective and projective \( C^*\)-cross norms on \( A \otimes B \) are the same for every \( C^*\)-algebra \( B \). (For separable \( C^* \)-algebras, this is equivalent to the property of being isomorphic to a \( C^* \)-subalgebra \( B \) of the Cuntz algebra \( \mathcal{O}_2 \) with the property that there exists a conditional expectation from \( \mathcal{O}_2 \) to \( B \).) However, the following result shows that, in contrast to the case of AF algebras, a solution to the classification problem for such \( C^* \)-algebras is necessarily less satisfactory.

\begin{theorem}[{\cite{Farah2014}}] \label{thm:nuclearturbulent}
The isomorphism relation  for unital simple separable nuclear
 \( C^* \)-algebras is not classifiable by countable structures, and is indeed strictly \( \leq_B \)-above the \( S_\infty \)-equivalence relation. 
\end{theorem}

In particular, simple nuclear \( C^* \)-algebras cannot be classified using only \( K \)-theoretic data. 
A similar anti-classification result has been obtained by Gardella and Lupini~\cite{Gardella2016PAMS} when considering e.g.\ the class of separable non-self-adjoint UHF operator algebras. (A \markdef{uniformly hyperfinite} (briefly, \markdef{UHF}) \markdef{algebra} is a \( C^* \)-algebra  which  is the direct limit \( \displaystyle{A = \lim_{\longrightarrow} (A_n , \varphi_n)} \) of an inductive system where each \( A_n \) is a 
finite-dimensional full matrix algebra and each \( \upvarphi_n \colon A_n \to A_{n+1} \) is a unital embedding.)

In~\cite{Farah2014} it was also shown that the isomorphism relation on unital simple separable nuclear
 \( C^* \)-algebras is   Borel reducible to the action of the automorphism group of the Cuntz algebra \( \mathcal{O}_2 \) on its closed subsets, and thus it is an essentially orbit equivalence relation. In~\cite{Elliot2013} the same upper bound was shown to hold for arbitrary separable \( C^* \)-algebras, as well as for many other related objects such as separable operator spaces or separable operator systems.

The search for the exact complexity of the classification problem for separable \( C^* \)-algebras culminated in the following result of Sabok, which in particular improves Theorem~\ref{thm:nuclearturbulent}. Recall that an \markdef{approximately interval} (or \markdef{AI}) \markdef{\( C^* \)-algebra} is a direct limit
\( \displaystyle{A = \lim_{\longrightarrow} (A_n , \varphi_n)} \) 
where, for each \( n \in \mathbb{N} \), \( A_n \cong F_n \otimes C([0;1]) \) for some 
finite-dimensional \( C^* \)-algebra \( F_n \) and \( \varphi_n \colon A_n \to A_{n+1} \) is a \( * \)-homomorphism.

\begin{theorem}[{\cite{Sabok2016}}] \label{thm:mainsabok}
The isomorphism relation on separable (nuclear, simple) \( C^* \)-algebras is complete for orbit equivalence relations. The same is true when restricting to the class of simple AI  \( C^* \)-algebras.
\end{theorem}

Interestingly enough, the proof of the above result goes through an analogous result for Choquet simplices. A few years earlier, Farah, Toms and T\"ornquist~\cite{Farah2014} showed
that the relation of affine homeomorphism on Choquet simplices is Borel
reducible to the isomorphism relation on simple, separable, AI \( C^* \)-algebras. In the same paper, they also showed that such relation is turbulent, yet Borel reducible to the orbit equivalence relation of a Polish group action.
Sabok's main contribution in~\cite{Sabok2016} was that of showing that the relation of affine homeomorphism of Choquet simplices
is indeed complete for orbit equivalence relations, yielding Theorem~\ref{thm:mainsabok} as a corollary.

Other subclasses of \( C^* \)-algebras on which isomorphism has maximal complexity are the following (this is a consequence of Theorem~\ref{thm:zielinski}):

\begin{theorem}[{\cite{Zielinski2016}}] \label{thm:zielinskiC^*algebras}
The isomorphism relations on \( C(K) \)-spaces and on commutative
\( C^* \)-algebras are both complete for orbit equivalence relations.
\end{theorem}

As for bi-embeddability between \( C^* \)-algebra we have:

\begin{theorem}[{\cite{Farah2014}}]
All \( \sigma \)-compact equivalence relations%
\footnote{An equivalence relation on a Polish space \( X \) is \(\sigma\)-compact if it can be written as a countable union of compact subsets of \( X^2 \).}
 are Borel reducible to bi-embeddability of unital AF \( C^* \)-algebras.
\end{theorem}

Since there are \( \sigma \)-compact equivalence relations which are not Borel reducible to any orbit equivalence relation, this implies that the relation of bi-embeddability on AF \( C^* \)-algebras is not an essentially orbit equivalence relation.

\subsubsection{Automorphism of \( C^* \)-algebras} \label{subsubsection:automorphisms}

If \( A \) is a separable \( C^* \)-algebra, the group \( \mathrm{Aut}(A) \) of \markdef{automorphisms of \( A \)} is a Polish group with respect to the topology of pointwise norm convergence. An automorphism of \( A \) is called (multiplier) inner if it is induced by the action by conjugation of a unitary element of the multiplier algebra \( M(A) \) of \( A \). 
Two automorphisms \( \alpha \) and \( \beta \) of \( A \) are called \markdef{unitarily equivalent} if \( \alpha \circ \beta^{-1} \) is inner or, equivalently, if for some unitary multiplier \( u \) one has
\( \alpha(x)=\beta(u x u^*) \) for all \( x \in A \). This defines a Borel equivalence relation on \( \mathrm{Aut}(A) \).

In the Eighties it was established that the automorphisms of 
\( A \) are concretely classifiable up to unitary equivalence  if and only if \( A \) has a \markdef{continuous trace}. (A separable \( C^* \)-algebra \( A \) has continuous trace if it is generated by abelian elements and the spectrum \( \widehat{A} \) is a Hausdorff space when endowed with the hull-kernel topology.) This equivalence has recently been turned into a subtle dichotomy by Lupini.

\begin{theorem}[{\cite{Lupini2014}}]
If \( A \) is a separable \( C^* \)-algebra that does not have continuous trace, then the automorphisms of \( A \) are not classifiable by countable structures up to unitary equivalence.
\end{theorem}

Thus either the classification problem for automorphisms of \( A \) up to unitary equivalence is smooth, or else it is very complicated and does not admit reasonably simple complete invariants.

Many results along the same lines have been proved in the subsequent years, for example:
\begin{itemizenew}
\item
If \( A \) is either the Jiang-Su algebra \( \mathcal{Z} \), or one of the Cuntz algebras \( \mathcal{O}_2 \) or \( \mathcal{O}_\infty \),  or a UHF algebra of infinite type, or a tensor product of a UHF algebra of infinite type and \( \mathcal{O}_\infty \), then the action of \( \mathrm{Aut}(A) \) on itself by conjugation is generically turbulent, i.e.\ all its orbits are meager and there exists a dense and turbulent orbit. It follows that if \( A \) is either \( \mathcal{Z} \)-stable (i.e.\ \( \mathcal{Z} \otimes A \cong A \)), a separable stable \( C^* \)-algebra, or a separable McDuff \( \mathrm{II}_1 \) factor,  then the relation of conjugacy on the set  of approximately inner automorphisms of \( A \) is not classifiable by countable structures. (See~\cite{Kerr2015}.)
\item
The automorphisms of the Cuntz algebra \( \mathcal{O}_2 \) are not classifiable by countable structures up to conjugacy or cocycle conjugacy. The same conclusion holds even when considering only automorphisms of \( \mathcal{O}_2 \) of finite order.
(See~\cite{Gardella2016}.)
\end{itemizenew}

\subsection{Ergodic measure preserving transformations} \label{subsec:measurepreserving}

Let \( (X,\mu ) \) be a standard Borel probability space and denote by \( \mathrm{Aut}(X,\mu) \) the group of all \markdef{measure preserving transformations} on \( X \) (endowed with the weak topology, which turns it into a Polish group). Given a countable group \( G \), the space \( \mathcal{A}(G,X,\mu) \) of \markdef{measure-preserving \( G \)-actions} on \( X \) is the set of all homomorphism from \( G \) to \( \mathrm{Aut}(X,\mu) \), which is Polish as a subspace of \( \mathrm{Aut}(X,\mu)^G \). 
The group \( \mathrm{Aut}(X,\mu) \) naturally acts on \( \mathcal{A}(G,X,\mu) \) (and thus, in particular, on the ergodic actions) by conjugation: two \( G \)-actions \( \sigma,\tau \in \mathcal{A}(G,X,\mu) \) are equivalent if and only if there is \( T \in \mathrm{Aut}(X,\mu) \) such that on a measure-one set the condition \( T(\sigma(g)(x)) = \tau(g)(T(x)) \) holds for all \( g \in G \). 

In~\cite{Foreman2004} it was proved that if \( G \) is amenable, then the ergodic measure-preserving actions of \( G \) on the unit interval \( [0;1] \) (equipped with the Lebesgue measure) cannot be classified by countable structures up to conjugation by elements of \( \mathrm{Aut}([0;1]) \). An important special case of this phenomenon is that the ergodic measure preserving transformations of \( [0;1] \) cannot be classified by countable structures up to conjugation via a measure preserving transformation (this is nicely complemented by Theorem~\ref{thm:FRW}).

Using an important result from~\cite{Kechris2001},
T\"ornquist showed that Foreman-Weiss' result can be further improved as follows.

\begin{theorem}[{\cite{Toernquist2009}}] 
Let \( G \) be a countable discrete group and let \( H \leq G \) be an infinite
subgroup. For any standard Borel probability space \( (X, \mu) \)  we have the following.
\begin{enumerate-(1)}
\item 
The measure preserving \( H \)-actions on \( (X, \mu) \) that can be extended to an ergodic
measure preserving \( G \)-action on \( (X, \mu) \) cannot be classified up to conjugacy by
countable structures.
\item 
If \( H \) either contains an infinite abelian subgroup or is normal in \( G \), then the
measure preserving \emph{ergodic} \( H \)-actions on \( (X, \mu) \) that can be extended to an ergodic
measure preserving \( G \)-action on \( (X, \mu) \) are not classifiable by countable structures up to conjugacy.
\end{enumerate-(1)}
\end{theorem}

\section{Classifying uncountable structures and non-separable spaces} \label{sec:uncountable}

The astonishing achievements of Borel reducibility in settling classification problems for countable structures and separable spaces naturally led to the idea of developing a similar framework to deal with uncountable structures and non-separable spaces. 
Basically, one fixes an infinite cardinal \( \kappa \) and then equips the space \( \pre{\kappa}{2} \) of all binary \( \kappa \)-sequences with the topology generated by the sets of the form
\[ 
\Nbhd^\kappa_s = \{ x \in \pre{\kappa}{2} \mid s \subseteq x \}
 \] 
for \( s \) a binary sequence of length \( < \kappa \). (Notice that when \( \kappa = \omega \) this is just the usual Cantor space.) All other topological notions are then updated accordingly: for example, one considers \( \kappa^+ \)-Borel sets (i.e.\ sets in the smallest \( \kappa^+ \)-algebra generated by open sets) and \( \kappa \)-analytic sets (i.e.\ continuous images of \( \kappa^+ \)-Borel subsets of \( \pre{\kappa}{2} \)) instead of the classical Borel and analytic sets. We refer the reader to the monographs~\cite{Friedman:2011nx, AM} for more on the basics of this theory, which is often called \emph{generalized descriptive set theory}.

On the one hand, the theory of the generalized Cantor space \( \pre{\kappa}{2} \) is quite different when moving to uncountable \( \kappa \)'s: even very basic results from descriptive set theory become either false or independent from \( \mathsf{ZFC} \). On the other hand, it is still possible to naturally code structures of size \( \kappa \) and spaces of density \( \kappa \) as elements of \( \pre{\kappa}{2} \) and then analyze the complexity of the related classification problems, which are \( \kappa \)-analytic equivalence relations on (subsets of) \( \pre{\kappa}{2} \), through \( \kappa^+ \)-Borel reducibility, i.e.\ the preorder \( \leq^\kappa_B \) induced by \( \kappa^+ \)-Borel measurable reductions. 

In the rest of this section we present a (non-exhaustive) sample of results that have been obtained so far within this framework.

\subsection{Model theory}

We consider uncountable models of countable complete first-order theories \( T \) and denote by \( \cong^\kappa_T \) the isomorphism relation on (codes for) \( \kappa \)-sized models of \( T \) . Shelah's celebrated Main Gap Theorem states that for such a theory \( T \), if \( T \) is not classifiable (i.e.\ \( T \) is unsuperstable, NDOP, and NOTOP) or if \( T \) is classifiable and deep, then it has \( 2^\kappa \)-many models of size \( \kappa \) up to isomorphism; if instead \( T \) is classifiable shallow of depth \( \alpha \) (in which case \(\alpha\) is a countable ordinal), then it has at most \( \beth_\alpha\left(|\gamma|^{2^{\aleph_0}}\right) \)-many \( \kappa \)-sized models up to isomorphism, where \( \gamma \) is such that \( \kappa = \aleph_\gamma \). Notice that the upper bound obtained in the second case (which is \( < 2^\kappa \) for many cardinals \( \kappa \)) depends on both the depth \(\alpha\) of the theory and  the cardinal \( \kappa \).

In~\cite{Friedman:2011nx}, Friedman, Hyttinen and Kulikov uncovered a deep connection between the mentioned Shelah's dividing line and the complexity of \( \cong^\kappa_T \), namely, \( T \) is classifiable shallow if and only if \( \cong^\kappa_T \) is \( \kappa^+ \)-Borel (for \( \kappa \) any cardinal satisfying \( \kappa^{< \kappa} = \kappa > 2^{\aleph_0} \)). Since \( \kappa^+ \)-Borel sets can naturally be arranged in a complexity hierarchy of length \( \kappa^+ \), this raises the question of whether the depth of a classifiable shallow theory \( T \) and the \( \kappa^+ \)-Borel rank of \( \cong^\kappa_T \) (i.e.\ the minimal ordinal \( \beta < \kappa^+ \) such that \( \cong^\kappa_T \) belongs to the \( \beta \)-level of the \( \kappa^+ \)-Borel hierarchy) are somewhat related. This led to the following descriptive set theoretic version of Shelah's Main Gap Theorem.

\begin{theorem}[{\cite{Mangraviti2018}}] \label{thm:mangra}
Let \( T \) be a countable complete theory \( T \) and \( \kappa \) be a cardinal satisfying \( \kappa^{<\kappa} = \kappa > 2^{\aleph_0} \). If \( T \) is classifiable shallow of depth \( \alpha \), then \( \cong^\kappa_T \) has \( \kappa^+ \)-Borel rank \( \leq 4 \alpha \); otherwise \( \cong^\kappa_T \) is not \( \kappa^+ \)-Borel at all.
\end{theorem}

Since  the depth \(\alpha\) of a classifiable shallow theory is always a countable ordinal (and thus much smaller than \( \kappa^+ > \aleph_2 \), under the above assumptions), Theorem~\ref{thm:mangra} yields a strong dichotomy: either \( \cong^\kappa_T \) is very low in the \( \kappa^+ \)-Borel hierarchy (more precisely, at a countable level which depends only on the depth \(\alpha\) and not on the cardinal \( \kappa \)), or else it is not \( \kappa^+ \)-Borel. Moreover, one can show that \( T \) is \( \kappa \)-categorical if and only if \( \cong^\kappa_T \) is open, if and only if \( \cong^\kappa_T \) is clopen. The proof of the latter result passes through another interesting correspondence between model-theoretic and topological notions: indeed, the completeness of a theory \( T \) translates to a density condition for the isomorphism classes of its \( \kappa \)-sized models. 

Another side result which may be obtained from the above theorems is that, for suitable cardinals \( \kappa \), if \( T \) is a classifiable shallow theory and \( T' \) is not, then \( \cong^\kappa_T \) is strictly \( \leq^\kappa_B \)-below \( \cong^\kappa_{T'} \). In the same vein, the following result shows that there is a large gap (with respect to \( \kappa^+ \)-Borel reducibility) between the isomorphism relations over models of classifiable and unclassifiable theories.

\begin{theorem}[{\cite{Moreno2017}}]
The following is consistent with \( \mathsf{ZFC} \). If \( T , T' \) are countable complete first-order theories with \( T \) classifiable and \( T' \) not, then \( \cong^\kappa_T \) is strictly \( \leq^\kappa_B \)-below \( \cong^\kappa_{T'} \) for any \( \kappa = \kappa^{< \kappa} = \lambda^+ \) with \( 2^\lambda > 2^{\aleph_0} \) and \( \lambda^{< \lambda} = \lambda \). Indeed, the partial order \( (\pow(\kappa), {\subseteq}) \) can be embedded in the \( \leq^\kappa_B \)-structure of equivalence relations strictly between \( \cong^\kappa_T \)  and \( \cong^\kappa_{T'} \).
\end{theorem}

Along the same lines  one gets:

\begin{theorem}[{\cite{Hyttinen2017}}]
Let \( T,T' \) be countable complete first-order theories.
Under suitable assumptions on \( \kappa \), if \( T \)  is classifiable and \( T' \) is stable with the orthogonal chain property (briefly, OCP), then \( \cong^\kappa_T \) is \( \kappa^+ \)-Borel reducible to \( \cong^\kappa_{T'} \). If in addition \( \mathsf{V=L} \), then \( \cong^\kappa_{T'} \) is even complete for all \( \kappa \)-analytic equivalence relations.
\end{theorem}

 The OCP implies that \(T' \) is unsuperstable and it is common among stable unsuperstable theories. For example, it is easy to find complete theories of abelian groups (or more generally elementary theories of ultrametric spaces) that have such property. 

Notice however that the tight correspondence between simplicity in terms of stability theory and \( \kappa^+ \)-Borel reducibility breaks down in certain models of set theory. For example

\begin{theorem}[{\cite{Hyttinen:2012fj}}] \label{thm:completekappa}
There is a theory \( T_{\omega+\omega} \) which is stable, NDOP and NOTOP, yet if \(\mathsf{V=L} \) the isomorphism relation \( \cong^\kappa_{T_{\omega+\omega}} \) is complete for all \( \kappa \)-analytic equivalence relations (for suitable cardinals \( \kappa \)).
\end{theorem}

However, in certain forcing extensions \( \cong^\kappa_{T_{\omega+\omega}} \)  is not even \( S_\kappa \)-complete, so the above failure may be more a defect of the universe \( \mathsf{L} \) rather than of the theory itself. 

Another instance of the fact that the complexity of \( \cong^\kappa_T \) may depend on the model of set theory we are working in is given by the complete theory  of dense linear orders without end-points \( DLO \). Indeed, \( \cong^\kappa_{DLO} \) is \( S_\kappa \)-complete if \( \kappa^{< \kappa} = \kappa \), but under \( \mathsf{V = L } \) it becomes even complete for all \( \kappa \)-analytic equivalence relations whenever \( \kappa \) is a successor of a regular cardinal~\cite{Hyttinen:2012fj}.

\subsection{Bi-embeddability between graphs and groups}

The problem of classifying uncountable structures up to bi-embeddability has been considered as well, although in this case the results are mostly on the side of anti-classification theorems.

\begin{theorem}[{\cite{Motto-Ros:2011qc} for weakly compact cardinals, \cite{Mildenberger2019} for the general case}]
Let \( \kappa \) be an uncountable cardinal satisfying \( \kappa^{< \kappa} = \kappa > \omega \). Then the bi-embeddability relation on \( \kappa \)-sized graphs is complete for all \( \kappa \)-analytic equivalence relations. The same is true for \( \kappa \)-sized (non-abelian) groups.
\end{theorem}

\begin{theorem}[{\cite{Calderoni2018g}}] \label{thm:calduncountable}
Let \( \kappa \) be an uncountable cardinal satisfying \( \kappa^{< \kappa} = \kappa > \omega \). Then the relation of bi-embeddability on torsion-free abelian groups of size \( \kappa \) is complete for all \( \kappa \)-analytic equivalence relations.
\end{theorem}

It is worth noticing that the techniques used to prove Theorem~\ref{thm:calduncountable} are completely different from those used to prove Theorem~\ref{thm:TFAcountable}.

\subsection{Metric structures and Banach spaces}

Finally, we mention some results concerning the complexity of classification problems for non-separable spaces.

\begin{theorem}[{\cite{AM,MottoRos2017h}}]
Let \( \kappa \) be any infinite cardinal.
\begin{enumerate-(1)}
\item
Complete ultrametric spaces of density character \( \kappa \) are classifiable up to isometry by \( \kappa \)-sized structures (that is, the isometry relation on such spaces is \( \kappa^+ \)-Borel reducible to isomorphism on \( \kappa \)-sized graphs). 
\item
The relation of isometry on discrete complete metric spaces of density character \( \kappa \) is \( S_\kappa \)-complete.
\end{enumerate-(1)}
\end{theorem}

Combining this with e.g.\ Theorem~\ref{thm:completekappa} one gets that if  \( \mathsf{V=L} \) then, for suitable cardinals \( \kappa \), the isometry relation on complete (locally compact) metric spaces of density character \( \kappa \) is complete for all \( \kappa \)-analytic equivalence relations. In contrast, assuming \( \mathsf{AD}+ \mathsf{V=L(\RR)} \) one gets that for every uncountable \( \kappa \) there is a clopen
equivalence relation \( E \) on \( \pre{\kappa}{2} \) such that there is no reduction (of any kind) of \( E \) to isometry between complete
ultrametric and/or discrete metric spaces of density character \( \kappa \) (see~\cite{MottoRos2017h}).

As for isometric bi-embeddability we have:
\begin{theorem}[{\cite{AM}}]
Let \( \kappa \) be an uncountable cardinal satisfying \( \kappa^{< \kappa} = \kappa > \omega \).
Then the relation of isometric bi-embeddability between complete metric spaces of density character \( \kappa \) is complete for all \( \kappa \)-analytic equivalence relations. Indeed, the same is true already when restricting to discrete spaces. 
\end{theorem}

Finally, similar results may be obtained for Banach spaces.

\begin{theorem}[{\cite{AM}}]
Let \( \kappa \) be an uncountable cardinal.
\begin{enumerate-(1)}
\item
The isomorphism relation for \( \kappa \)-sized graphs is \( \kappa^+ \)-Borel reducible to linear isometry on Banach spaces with density \( \kappa \). Thus the latter relation is consistently complete for all \( \kappa \)-analytic equivalence relations (for suitable \( \kappa \)'s).
\item
If \( \kappa^{< \kappa} = \kappa > \omega \), then the relation of isometric bi-embeddability on Banach spaces with density \( \kappa \) is complete for all \( \kappa \)-analytic equivalence relations.
\end{enumerate-(1)}
\end{theorem}

\newcommand{\etalchar}[1]{$^{#1}$}

\end{document}